\documentclass[10pt, twocolumn, final]{IEEEtran}
\usepackage{overpic}
\usepackage{url}
\usepackage{subfigure}
\newtheorem{lemma}{Lemma}[section]

\include{epsf}
\usepackage{times}
\usepackage{latexsym}
\usepackage{amsmath}
\usepackage{amsfonts}
\usepackage{pstricks}
\usepackage{amssymb}
\usepackage{amsxtra}
\usepackage{subfigure}
\usepackage{epsfig}
\usepackage{color}
\usepackage{graphics}
\usepackage[ruled]{algorithm2e}

\usepackage{cite}
\DeclareMathOperator*{\argmin}{arg\,min}

\begin{document}
\title{Delay Estimation and Fast Iterative Scheduling Policies for LTE Uplink
\author{Akash Baid(Rutgers), Ritesh Madan(Qualcomm), and Ashwin Sampath(Qualcomm) }
}

\maketitle \psfull

\begin{abstract}

We consider the allocation of spectral and power resources to the mobiles (i.e., user equipment (UE)) in a cell every subframe (1 ms) for the Long Term Evolution (LTE)
orthogonal frequency division multiple access (OFDMA) cellular network. To enable scheduling based on packet delays, we design a novel mechanism for inferring the packet delays approximately from the buffer status reports (BSR) transmitted by the UEs; the BSR reports only contain queue length information. We then consider a constrained optimization problem with a concave objective function -- schedulers such as those based on utility maximization, maximum weight scheduling, and recent results on iterative scheduling for small queue/delay follow as special cases. In particular, the construction of the non-differentiable objective function based on packet delays is novel. We model constraints on bandwidth, peak transmit power at the UE, and the transmit power spectral density (PSD) at the UE due to fractional power control. When frequency diversity doesn't exist or is not exploited at a fast time-scale, we use subgradient analysis to construct an $O(N \log L)$ (per iteration with small number of iterations) algorithm to compute the optimal resource allocation for $N$ users and $L$ points of non-differentiability in the objective function. For a frequency diversity scheduler with $M$ sub-bands, the corresponding complexity per iteration is essentially $O(N(M^2+L^2))$. Unlike previous iterative policies based on delay/queue, in our approach the complexity of scheduling can be reduced when the coherence bandwidth is larger. Through detailed system simulations (based on NGMN and 3GPP evaluation methodology) which model H-ARQ, finite resource grants per sub-frame, deployment, realistic traffic, power limitations, interference, and channel fading, we demonstrate the effectiveness of our schemes for LTE.
\end{abstract}

\section{Introduction}

Wideband cellular systems such as LTE allow for resource allocation with high granularity of a resource block (RB) of 1 ms by 180 KHz~\cite{lte_book}. While control signalling and the general framework for the physical and medium access control (MAC) layers is specified to enable efficient use of spectral resources, the exact resource allocation algorithms for power and frequency allocation can be designed by an implementor. Moreover, each cell can serve on the order of a thousand active connections over a bandwidth of 20~MHz. Hence, in order to take advantage of the flexibility allowed in resource allocation, the resource allocation algorithms have to be computationally simple. Many schedulers in the literature entail maximizing the weighted sum of rates in each subframe. For example, the weights could be based on utility functions of average rate~\cite{kushner_2004},\cite{stolyar_2005}, the queue length~\cite{andrews_2004},~\cite{tassiulas_1990}, or head-of-line delay~\cite{shakkottai_2002},~\cite{sadiq_2010}. In the uplink, the resource allocation problem must consider the maximum transmission power of a mobile and the constraints on the transmission power imposed by fractional power control to limit inter-cell interference~\cite{lte_book},~\cite{fpp_lte}. When contiguous bandwidth allocation is considered, the problem of maximizing the weighted sum rate in each subframe on the UL can be posed as a constrained convex optimization problem. For $N$ users and $M$ sub-bands general purpose methods can solve the problem in $O((NM)^3)$. With peak UE power constraints, a $O(NM)$ per iteration subgradient algorithm was obtained in~\cite{huang_2009}; heuristics to compute allocations with integral number of resource blocks (RBs) were considered as well. Interior point methods (which have faster convergence) with an $O(NM^2)$ (if $N>>M$) Newton iteration were obtained in~\cite{madan_2011} for uplink resource allocation with additional fractional power control constraints. However non-differentiable objective functions are not considered under the framework in~\cite{madan_2011}. 

Also relevant to our paper are recent results on low complexity iterative scheduling algorithms. Many papers prior to these results had considered scheduling to maximize the  sum of weighted rates in subframe $n$, where the weights were based on the arrivals and departures in the queue of a user until subframe $n-1$. The iterative policies in~\cite{sharma_2011,bodas_2011} take into account how the weights change in subframe $n$ to determine the resource allocation in \emph{that} subframe. In particular, the queue based server side greedy (SSG) rule is proposed for multi-rate channels in~\cite{bodas_2011} and a delay based rule with iterative matching in each subframe for ON-OFF channels is considered in~\cite{sharma_2011}. The results in these papers shed a remarkable insight that when the rate grows linearly with bandwidth (no peak power constraints at the transmitter), as the number of users in the system grow, these rules lead to much smaller per-user queues and delays, respectively, compared with previous approaches. However, the complexity of these algorithms grow with the resource granularity even if the coherence bandwidth does not grow. In this paper, we construct a continuous but non-differentiable concave reward function based on packet delays. We argue that the matching algorithm in~\cite{sharma_2011} is an approximate algorithm to maximize this reward function in every subframe.

Motivated by the above observation, we consider resource allocation to maximize a continuous (possibly non-differentiable) concave reward function. We first consider a channel model where the channel gain in the frequency domain is flat and formulate the resource allocation problem as a non-differentiable convex optimization problem. Note that in typical cellular environments, the channel gains can be fairly correlated even for frequencies 2 to 5 MHz apart~\cite{zhang_2007} -- hence, the assumption of frequency flat fading is a reasonable one when the total bandwidth is up to 5 MHz (28 RBs) or lower, or if the UEs are allocated to sub-bands ($<$ 5 MHz) over a slower time-scale based on interference and channel statistics.
The above assumption allows us to use subgradient analysis to design algorithms with $O(N\log L)$ cost per iteration (with small number of iterations) for $N$ users and $L$ points of non-differentiability in the objective function.
We discuss implementation issues for the resulting algorithm in a practical LTE system with H-ARQ re-transmissions, finite number of resource grants per subframe, and the constraint that all uplink transmissions have to be over a contiguous set of RBs. Notably, we also design a novel mechanism to estimate head-of-line delays of queues at UEs with low complexity via only queue length information contained in the buffer status reports (BSR). We note our techniques are equally applicable for enabling delay based scheduling in the PCF and HCF modes in WiFi~\cite{spec_802_11}. We demonstrate the improvement in performance due to our techniques through numerical results obtained via comprehensive numerical simulations based on 3GPP evaluation methodology~\cite{eval_tr}. Finally, when frequency selective fading is considered, we show how interior point methods with complexity of $O(NM^2 + NL^2)$ per Newton iteration can be obtained; note that in practice $N>>L,M$. Since, we consider non-differentiable cost functions, this requires additional analysis compared to that in~\cite{madan_2011} where only differentiable cost/utility functions were considered.

Prior work in devising practical resource allocation schemes for the LTE uplink includes:
\cite{calabrese_2007} considers allocation of fixed size resource chunks to UEs, \cite{calabrese_2008b} extends this
approach where each (RB,UE)-tuple is associated with a metric (which cannot capture power constraint at power limited
UE), a similar (RB,UE) metric is considered in~\cite{lee_2009}. These methods do not extend to solving a general
resource allocation problem considered in this paper.
Semi-persistent scheduling for voice over IP (VoIP) has been considered in, for example,~\cite{jiang_2007}.
In~\cite{yaacoub_2009} heuristics for maximizing utilities of UEs in each subframe in the presence of frequency selective fading but no fractional power control were considered. Similarly, heuristics to satisfy minimum rate constraints of most users and maximize the sum rate were considered in~\cite{gao_2008}, heuristics to maximize sum weighted rate were designed in~\cite{kim_2005}, and algorithms for long term proportional fairness were considered in~\cite{ma_2009}.

\section{System Model}
\label{sec:system}
\subsection{Channel Model, Power, Rate}
\label{subsec:frac_power}
We focus on the uplink of a single cell in LTE with $N$ UEs and the total bandwidth divided
into $M$ sub-bands of equal bandwidth $B$, with $B$ less than the coherence bandwidth of each user.
The maximum transmit power of each UE is $P$.
The channel gain for
UE $i$ on sub-band $j$ is $G_{ij}$; we focus on the
scheduler computation in a subframe, and don't explicitly show
the dependence of quantities on time $t$. The base-station can measure the $G_{ij}$s via
decoding the sounding reference signal (SRS)~\cite{lte_book}. Fractional power control in LTE limits the amount of interference
a UE causes at base-stations in neighboring cells. A UE which is closer to the cell edge
inverts a smaller fraction of the path loss to the serving base-station than a UE
which is closer to the serving base-station~\cite{fpp_lte}. Thus the transmit powers of a UE on different
sub-bands satisfy~\cite{madan_2011}:
\[
\smash{{p_{ij} \leq \gamma_{ij} b_{ij}, \forall i,j, \quad \sum_{j=1}^M p_{ij} \leq P,}}
\]
where $b_{ij}$ is the bandwidth allocated to UE $i$ on sub-band $j$ and $\gamma_{ij}$ is a sub-band specific constant.

The interference PSD at the serving base-station
on sub-band $j$ (denoted as $I_j$)  can be
measured by the base-station periodically over unassigned frequency
resources. The value depends on the
interference coordination algorithm used~\cite{fodor_2009}.
When a UE transmits with power $p_{ij}$ over bandwidth $b_{ij}$ on
sub-band $j$, it achieves a rate given by (treating interference as noise)
\[
\smash{b_{ij}\psi \left( \frac{G_{ij}p_{ij}}{ b_{ij}I_j}\right)}
\]
where $\psi: \mathbb{R}_+\mapsto \mathbb{R}_+$ is an increasing concave and differentiable function which maps the SINR to spectral efficiency.


\subsection{Control Signaling}
\label{subsec:bsr_sr}
Single carrier frequency division multiple access (SC-FDMA) is used in the LTE uplink~\cite{lte_book} and so a UE can be granted a number of 180 kHz resource blocks in a contiguous manner in frequency. The resource allocation to the UEs is computed by the base-station every subframe (1~ms) and signalled to the UEs via resource grants which include the contiguous set of RBs allocated to the UE and the modulation and coding scheme (MCS). The timeline is as follows: a resource grant is made to the UE at time $t$ for an uplink
transmission at time $(t+4)$. At time $(t+8)$ the base-station transmits an ACK/NACK to indicate if it could decode the packet; if a NACK is received by the UE, it re-transmits at the same power and uses the same RBs at time $(t+12)$ (as
at time $(t+4)$). We assume a constant number of maximum allowable re-transmissions for all UEs and do not
adapt the re-transmission power and resource assignment through additional control signalling available in LTE.

Buffer status report (BSR) and scheduling request (SR) are transmitted by the UEs to inform the base-station about new packet arrivals at the UE. We
describe the mechanism for the special case of single logical channel (LC), or connection, at each UE.
SR is one bit of information used to indicate the arrival of packets in an empty buffer at the UE.
Each UE periodically
gets an opportunity to send SR, and the time interval between
two successive opportunities for SR is denoted by $T^\text{SR}$, and is assumed to be fixed
in a cell.
BSRs contain a quantized value of the number of bytes pending transmission at the UE\footnote{We ignore
the effect of quantization in BSR, but the methods in this paper extend easily to quantized BSR.}, and are generated in two different ways:
\emph{Regular BSR:} If the queue is empty in subframe $t$, and new packets arrive
in subframe $t+1$, a regular BSR is generated at time $t+1$. When a regular BSR is generated, a SR is transmitted at the next available SR opportunity unless resources are granted to the UE between the BSR generation and the opportunity to transmit SR.
\emph{Periodic BSR:} A periodic BSR is generated every $T^\text{BSR}$ subframes.
A periodic BSR thus generated is transmitted by the UE to the base-station at the earliest subframe after generation when resources are granted to it by the base-station.

\vspace{-0.3cm}
\begin{figure}[!th]
\par
\begin{center}
\includegraphics[width=0.30\textwidth]{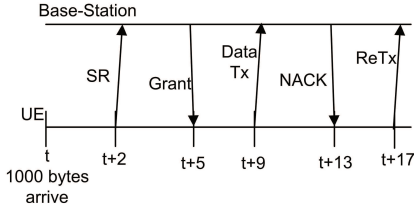}
\caption{Uplink timeline}
\label{fig:ulframe}
\end{center}
\end{figure}
\vspace{-0.3cm}

A typical sequence of transmissions is shown in Fig.~\ref{fig:ulframe}. The buffer at a UE is empty in subframe
$(t-1)$ and a 1000 byte packet arrives in subframe $t$. The next SR opportunity is subframe $(t+2)$ -- the SR transmission
by the UE signals to the base-station that the buffer at the UE is non-empty. In response, the base-station allocates resource
to the UE on the uplink via a grant at time $(t+5)$ -- the actual uplink transmission occurs 4 subframes later, i.e., in subframe
$(t+9)$. This transmission includes the BSR report. Assume that the UE is allocated enough resources to also transmit 200 bytes of the data packet -- then the BSR report will contain a value of $800$ bytes for the left-over data at the UE. The first transmission is unsuccessful -- this is indicated by a NACK transmitted by the base-station at time $(t+13)$. The UE re-transmits the packet at time $(t+17)$ -- this transmission is decoded successfully by the base-station, and hence it is known at the base-station that 800 bytes were pending transmission at the UE at time $(t+9)$ which
is the time when the BSR report was created.

\section{Reward Functions }
\label{sec:reward}
In this section, we define the reward functions that we use for the optimization problem and relate it to the schemes used in earlier works. We assume each UE to have one active LC which supports either best effort or delay QoS traffic.\\

\subsection{Best Effort}
A flow, $i$, which is
best-effort is associated with an average rate
\mbox{$x_i(t)\in \mathbb{R}_{+}$} in subframe $t$ which is updated as follows:
\begin{equation}
\label{eqn:be_avg}
\smash{x_i(t+1) = (1-\alpha_i)   x_i(t) + \alpha_i r_i, \quad \forall t\geq 0,}
\end{equation}
where $r_i$ is the rate at which UE $i$ is served in the current subframe, and $0<\alpha_i < 1$ is a user specific constant. The user
experience in subframe $t$ is modeled as a strictly concave increasing
function \mbox{$U_i:\mathbb{R}_+ \mapsto\mathbb{R}$} of the average
rate $x_i(t)$. Traffic for applications such as file transfer and web browsing can be modeled by best effort flows, and is typically transferred over a TCP connection which has closed loop rate control. We greedily maximize the total utility at each time-step, i.e., the reward function for UE $i$ with best
effort traffic, at time
$t$ is~\cite{madan_2010}
\begin{equation}
\label{eqn:be_rew}
f_i(r_i) = \frac{1}{\alpha_i}U_i((1-\alpha_i)x_i(t) + \alpha_i r_i).
\end{equation}
If we set $f_i(r_i) = U'(x_i(t))r_i$, and let $\alpha_i\rightarrow 0$ in equation~(\ref{eqn:be_avg}), the resulting scheduler is identical to that in~\cite{stolyar_2005}. Thus, our analysis offers a computationally efficient method to implement the scheduling policy in~\cite{stolyar_2005} for the LTE uplink with fractional power control; also note that it is easy to show that the rate vectors in the uplink resource allocation problem satisfy the conditions required for the results in~\cite{stolyar_2005}.

\subsection{Delay QoS Traffic}
\label{subsec:rew_delay}
 Here the user experience is a function of the packet delays. User experience is acceptable when the packet delays are lower than a certain tolerable value. The packet arrival process
is assumed to be independent of the times at which the packets are served. Traffic for applications such as voice calls and
live video chatting fall in this category.

At time $t$, let $\pi_i(t)$ be the number of packets in the queue of UE $i$.
Denote the sizes and the delays of these $\pi_i(t)$ packets by $\{s_i(1), \hdots, s_i(\pi_i(t))\}$ and
$\{d_i(1), \hdots, d_i(\pi_i(t))\}$.
Then for a UE $i$ with delay QoS traffic, we define the reward function as:
\begin{equation}
\label{eqn:rew_qos}
\begin{aligned}
f_i(r_i) &= \sum_{j=1}^{n_i^\text{serv}(r_i)} s_i(j)d_i(j)\\
& + \left(r_i\Delta - \sum_{j=1}^{n_i^\text{serv}(r_i)}s_i(j)\right)d_i(n_i^\text{serv}(r_i)+1)
\end{aligned}
\end{equation}
where $\Delta$ is the length of a subframe (1 ms) and $n_i^\text{serv}(r_i)$ is the number of packets from UE $i$ served fully if UE $i$ is scheduled at rate $r_i$, i.e., \\
$n_i^\text{serv}(r_i) = \max\left\{ k: \sum_{j=1}^k s_i(j) \leq r_i\Delta\right\}$.

\begin{lemma}
$f_i(r_i)$ is a continuous concave function.
\end{lemma}
\begin{proof}
Concavity follows from the observation that $d_i(1)> \hdots> d_i(\pi_i(t))$ and continuity is immediate
from definition.
\end{proof}

\vspace{-0.3cm}
\begin{figure}[!th]
\par
\begin{center}
\includegraphics[width=0.40\textwidth]{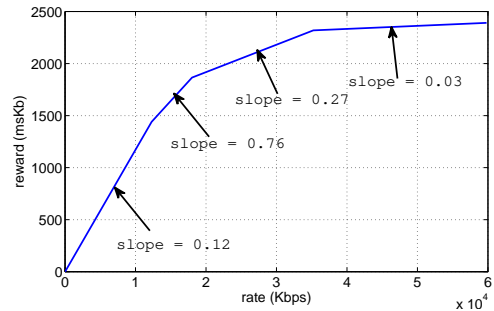}
\caption{Example reward function for delay QoS flow.}
\label{fig:reward_eg}
\end{center}
\end{figure}
\vspace{-0.3cm}

\emph{Example:} Consider a UE with delay QoS traffic and four packets in the queue with delays (in ms) at
time $t$ given by
$ d_1 = 120, d_2 = 76, d_3 = 27, d_4 = 3,$
and packet sizes (in KB) are
$ s_1 = 1.5, s_2 = 0.7, s_3 = 2.1, s_4 = 3$.
Then the corresponding reward function $f_i$ is shown in Fig.~\ref{fig:reward_eg}.

\subsection{Iterative Queue and Delay Based Policies}
If we restrict the model in~\cite{sharma_2011} to frequency flat fading, i.e., a user is either connected to no server or all servers at any time,
the algorithm in that paper can be interpreted as one which approximately maximizes the reward function in equation~(\ref{eqn:rew_qos}). Specifically, the matching algorithm reduces to one where in each iteration a server is allocated to a user with the highest head-of-line delay times spectral efficiency -- this approximately equalizes the head-of-line delay times spectral efficiency for all users after the allocation, which (as we will show) is the optimality condition to maximize the reward function in~(\ref{eqn:rew_qos}) when divisible servers are considered. Larger the number of servers the same bandwidth $B$ is divided into, the closer the approximation. Note that peak power constraints are not modeled in~\cite{sharma_2011}. When frequency selective fading is considered, i.e., a user may be connected to a subset of servers in the model in~\cite{sharma_2011}, there is a sequence of maximum weight matchings which will approximately compute a solution which maximizes the reward function in~(\ref{eqn:rew_qos}). Motivated by this interpretation, we consider the maximization of the reward function in~(\ref{eqn:rew_qos}) for a much more general model with multiple rate options, peak power constraints, and different transmit PSD constraints on different sub-bands. We also note that the complexity of the algorithm in~\cite{sharma_2011} is $O(NR^2)$ for $N$ users and $R$ RBs -- when there are multiple RBs in each sub-band of bandwidth $B$, the complexity of our algorithms is lower. Finally, similar connections can be drawn between the scheme in~\cite{bodas_2011} for frequency flat fading and using an objective function based on sums of squares of queues as in~\cite{sadiq_2009}; the connections for the frequency selective fading case seem to exist but are harder to analyze.

\section{Estimation of Packet Delays}
\label{subsec:inferring_arrivals}
We now describe a method to infer approximate packet delays at the eNB via the mechanisms available in LTE. We use the SR and BSR report generated at the start of the burst of packets, and periodic BSR reports which are generated regularly but transmitted only when resources are allocated to the UE (see Sec.~\ref{subsec:bsr_sr}), along with the scheduling decisions made by the base-station to estimate packet delays. The main intuition is as follows: if the base-station estimates the queue length at time $t$ to be say, 1000 bytes, but later decodes a BSR which was created at time $t$ and has value $1300$ bytes, the base-station can deduce that 300 bytes arrived between time $t$ and the time at which the previous BSR was created. This information about the time interval during which the 300 bytes arrived can be used for making resource allocation decisions -- specifically, scheduling policies based on packet delays can be implemented. The main complexity is due to re-transmissions which can lead to the BSR report arriving out of order at the base-station.

Let $T^\text{retx}$ be the maximum amount of time between the first transmission of a MAC packet and the latest time when it can be re-transmitted
for H-ARQ (for example, if we configure 6 as the maximum number of re-transmissions, $T^\text{retx} = 48$ subframes). We estimate the number of bytes that arrived, $A_i(t)$ in each subframe $t$.
The buffer status reports are denoted by a sequence of random three tuples:
$$\smash{\{B_i(1), \tau_i(1), \delta_i(1)\}, \{B_i(2), \tau_i(2), \delta_i(2)\}, \hdots}$$
where $B_i(1)$ is the buffer size reported in first BSR, $\tau_i(1)$ is the time at which first BSR was received, and $(\tau_i(1) - \delta_i(1))$ is the time at which
the first BSR was generated, and so on. $C_i(t)$ denotes the number of bytes scheduled for transmission from UE $i$, $\hat{C}_i(t)$ the
number of bytes which were successfully received from UE $i$, and $F_i(t)$ the number of bytes that failed the final re-transmission for UE $i$, at time $t$.

We maintain the history of estimated queue length for each UE $i$ for duration $T^\text{retx}$, denoted by $Q_i(t-T^\text{retx}:t)$. Then, we update
the $Q$ matrix and the arrival vector $A$, at each $t$ as follows:

For every $t$, $i$
\begin{enumerate}
\item \emph{Scheduled Bytes:} $Q_i(t) = Q_i(t-1) - C_i(t)$.
\item \emph{Failed Bytes:} $Q_i(t) = Q_i(t) + F_i(t)$.
\item \emph{BSR report:} If a BSR report is received at time $t$, i.e., there is $n$ such that $\tau_i(n) = t$, then update queue
state as follows:
If the base-station has not received any BSR report created after time $t-\delta_i(n)$, then
$$Q_i(t-\delta_i(n):t) = Q_i(t-\delta_i(n):t) + A_i(t-\delta_i(n))$$ where
arrival $A_i(t-\delta_i(n)) =  B_i(t) - Q_i(t-\delta_i(n))$
otherwise for
$$\argmin_{\{m:~\tau_i(m) < t\}}  [\tau_i(m) - \delta_i(m) - (\tau_i(n)-\delta_i(n))]$$
update
\[\begin{aligned}
&A_i(t-\delta_i(n)) =  B_i(t) - Q_i(t-\delta_i(n))\\
&A_i(\tau_i(m)-\delta_i(m)) = A_i(t-\delta_i(m)) - A_i(t-\delta_i(n))\\
\end{aligned}
\]
\[\begin{aligned}
&Q_i(t-\delta_i(n):\tau_i(m)-\delta_i(m) - 1)\\
& = Q_i(t-\delta_i(n):\tau_i(m)-\delta_i(m) - 1) + A_i(t-\delta_i(n))
\end{aligned}
\]
\end{enumerate}
Note that $Q_i$ can have negative entries.

\section{Frequency Flat Fading}
\label{sec:freq_flat_fading}

Here, we consider the resource allocation to $N$ UEs over a single sub-band with bandwidth $B$ and frequency flat fading. We drop the dependence of quantities in the general model on the sub-band $j$ -- for example, we denote channel gain from UE $i$ to the eNB as $G_i$. We allow for contiguous allocation -- this is a reasonable approximation when $B$ is larger than a few RBs. Rounding techniques in, for example,~\cite{huang_2009} can be used to obtain integral solutions. The optimization problem to maximize the sum of rewards for all UEs over the bandwidth allocation vector
$b \in\mathbb{R}_+^N$ in a subframe is:
\begin{equation}
\label{opt:cvx}
\begin{aligned}
\text{max.} \quad & \sum_{i=1}^N f_i\left(b_i\psi\left(\frac{G_i \min(\gamma_i b_i, P) }{Ib_i}\right) \right)\\
\text{s.t.} \quad & 0\leq b_i\leq b_i^\text{max}, ~~\forall i, \quad \sum_{i=1}^N b_i \leq B\\
\end{aligned}
\end{equation}
where $b_i^\text{max}$ is the maximum bandwidth that UE $i$ can use based on
the estimated queue length, $Q_i(t)$, for UE $i$, and satisfies:
\[
 b_i^\text{max}\psi\left(\frac{G_i \min(\gamma_i b_i^\text{max}, P)}{Ib_i^\text{max}}\right) = Q_i(t)/\Delta
\]
where we recall that $\Delta$ is the length of a subframe (1 ms). Since, the function on the left is an increasing function
of $b_i^\text{max}$, we can compute $b_i^\text{max}$ efficiently via a bisection search. Problem~(\ref{opt:cvx}) is a convex optimization problem (with non-differentiable objective function) due to the lemma which follows.
\begin{lemma}
\label{lem:cvx}
The objective function in optimization problem~(\ref{opt:cvx}) is concave in the $b_i$s for $b_i \geq 0$, for all $i$.
\end{lemma}
\begin{proof}
Consider the function $g:\mathbb{R}_+\mapsto \mathbb{R}_+$ defined by
$
g(x) = x\psi(c/x),~\forall x > 0, c \in \mathbb{R}_+ \text{  is constant}.
$
Since, $\psi$ is assumed to be concave, it is easy to verify (via showing that the second derivative is always negative)
that $g$ is concave as well.
Since, (i) the sum of concave functions is concave, and (ii) the composition of one concave function with another is concave,
to show that the objective function is concave, it is sufficient to show that the following function is concave
\[
h(x) = x \psi\left(\frac{\min(c_1 x, c_2)}{x}\right),~ \forall x \geq 0, c_1,c_2 \in \mathbb{R}_+ \text{ are constant}
\]
Note that the above function is well defined for $x\geq 0$. Since, $\psi$ is an increasing function, we can write
$
h(x) = \min \left\{ x \psi(c_1), x\psi(c_2/x) \right\},
$
which is the minimum of two concave functions, and hence, concave.
\end{proof}

\subsection{Characterization of Optimal Solution}
We define a function which maps the bandwidth allocation $b_i$ to achievable rate for user $i$:
$$
\smash{h_i (b_i) = b_i\psi\left(\frac{G_i \min(\gamma_i b_i, P) }{Ib_i}\right) }
$$
We denote the sub-differential of a function $g:\mathbb{R}\mapsto\mathbb{R}$ at $x$ by $\partial g(x)$. For continuous
concave functions over the set of reals, the subdifferential at $x$ is the set of slopes of lines tangent to $f$ at $x$.

Let $b^\star\in \mathbb{R}_+$ denote the solution to the resource allocation problem~(\ref{opt:cvx}).
The following lemma shows that an optimal allocation in a given subframe is one for which the following quantities are equal for all users with non-zero bandwidth allocation: for best effort user, the marginal utility times the incremental rate when more bandwidth is allocated to it, and for delay QoS user, the delay of the oldest packet which is not served completely times the incremental rate when more bandwidth is allocated to it.

\begin{lemma}
\label{lem:opt}
There
exists a $\lambda^\star>0$ such that if $i$ is best effort, then
$$\begin{aligned}
\lambda^\star \in  U'((1-\alpha)x_i(t) + \alpha_i r_i^\star)\partial h_i(b_i^\star),&\quad\text{if } b_i^\star > 0\\
\lambda^\star < U'((1-\alpha)x_i(t) )\min \partial h_i(0),~~\quad&\quad\text{if } b_i^\star = 0\\
\end{aligned}
$$
else, if $i$ is delay QoS and $b_i^\star>0$,
\begin{itemize}
\item if $\sum_{j=1}^{n_i^\text{serv}(r_i^\star)} s_i(j) < r_i^\star \Delta$, ~~
$\lambda^\star \in d_i(n^\text{serv}_i(r_i^\star)+1)\partial h_i(b_i^\star)$
\item else if $\sum_{j=1}^{n^\text{serv}(r_i^\star)} s_i(j) = r_i^\star \Delta$
$$\lambda^\star \in \left[d_i(n^\text{serv}_i)\min \partial h_i(b_i^\star),  d_i(n^\text{serv}_i+1)\max \partial h_i(b_i^\star)\right]$$
\end{itemize}
else, if $i$ is delay QoS and $b_i^\star=0$,
$$
\lambda^\star < d_i(1)\min \partial h_i(0)
$$
where $r_i^\star = h_i(b_i^\star)$.
\end{lemma}
\begin{proof}
The lemma follows from standard arguments in, for example~\cite{shor_1985}, the definitions of $f_i$'s, and that the subdifferential of $f_i$ for delay QoS user $i$ is given by
\[
\hspace{-0.1cm}\partial f_i(r_i)=\left \{
\begin{array}{ll}
d_i(n^\text{serv}_i+1),& \sum_{j=1}^{n_i^\text{serv}(r_i^\star)} s_i(j) < r_i^\star \Delta \\
\left[d_i(n^\text{serv}_i),  d_i(n^\text{serv}_i+1)\right],
& \sum_{j=1}^{n^\text{serv}(r_i^\star)} s_i(j) = r_i^\star \Delta
\end{array}\right .
\]
\end{proof}

We now evaluate the sub-differential of $h_i$ for $x\geq 0$, which is bounded because $\gamma_i$ is assumed to be bounded.
\[
\partial h_i (x) = \left \{
\begin{array} {l}
\left\{\psi \left(\frac{G_i(t) \gamma_i}{I}\right)\right \},~\text{if } x < P/\gamma_i\\
\\
\left\{\psi \left(\frac{G_i(t) P}{Ix}\right) -\frac{G_i(t)P}{x}\psi'\left(\frac{G_i(t) P}{Ix}\right)  \right \},\text{if } x > P/\gamma_i\\
\\
\left[\psi \left(\frac{G_i(t) \gamma_i}{I}\right)- \frac{G_i(t)P}{x}\psi'\left(\frac{G_i(t) \gamma_i}{I}\right),\right .
 \\\qquad\qquad \left.\psi\left(\frac{G_i(t) \gamma_i}{I}\right)\right ],~
\text{if } x = P/\gamma_i\\
\end{array}
\right .
\]

\vspace{-0.3cm}
\begin{figure}[!th]
\par
\begin{center}
\includegraphics[width=0.50\textwidth]{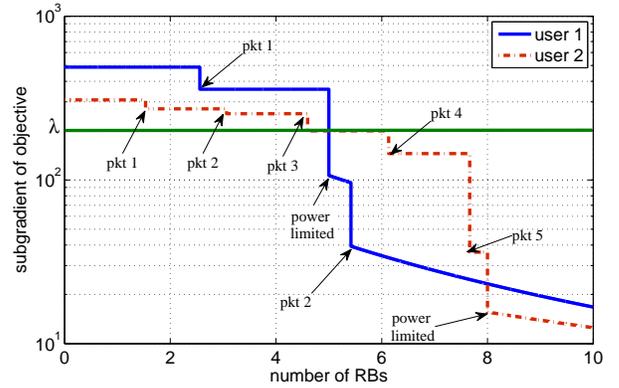}
\caption{Optimality condition}
\label{fig:opt_eg}
\end{center}
\end{figure}
\vspace{-0.3cm}

We illustrate the optimality condition via a two user example. The total bandwidth to be shared is 10 RBs, or 1800 KHz.
All packets are of size 500 bits. The packet delays of the two users in the given subframe are
\[
\begin{aligned}
\smash{\text{User 1:}}&\smash{\quad [450, 330, 135, 80, 20]}\\
\smash{\text{User 2:}}&\smash{\quad [170, 150, 140, 110, 80, 20]}\\
\end{aligned}
\]
The rate at which the users can be served as a function of the RBs are given by:
\[
\begin{aligned}
h_1(b_1)=&\left\{\begin{array}{ll}
                b_1\log_2\left(1+ 10^{0.05}\right) & b_1 \leq 5*180 \text{khz}\\
                b_1\log_2\left(1+ 10^{0.05}\frac{5*180}{b_1}\right) & b_1 > 5*180 \text{khz}\\
                  \end{array}
  \right.\\
h_2(b_2)=&\left\{  \begin{array}{ll}
                b_2\log_2\left(1+ 10^{0.4}\right) & b_2 \leq 8*180 \text{khz}\\
                b_2\log_2\left(1+ 10^{0.4}\frac{8*180}{b_2}\right) & b_2 > 8*180 \text{khz}\\
                  \end{array}
                  \right.
\end{aligned}
\]
where the 5 and 8 RB thresholds (and corresponding SINRs of 0.5~dB and 4~dB) are derived from fractional power control constraints in Section~\ref{subsec:frac_power}. The subgradient of the rewards for both the users as a function of bandwidth allocation, and the optimal bandwidth allocation are shown in Fig~\ref{fig:opt_eg} -- the optimal resource allocation is $5$ RBs to each user, and the optimal dual variable $\lambda^\star$ is shown in the figure. For each user, the figure also shows the number of RBs required to fully serve a given number of packets and the number of RBs at which the user becomes power limited, i.e., the maximum peak power constraint limits the transmission power rather than the fractional power control which limits the transmit PSD.


\subsection{Computation of Optimal Solution}
\label{subsec:bisect}
The optimization problem~(\ref{opt:cvx}) entails the maximization of the sum of concave functions subject to a linear inequality constraint. While, in principle, the optimal resource allocation scheme can be computed via a bisection search on
the dual variable $\lambda$, two difficulties arise: (i)~There may be multiple values of $b_i$ for which the subgradient of $f_i\circ h_i$ is equal to $\lambda$. See, for example, the first packet for user~1 in Fig.~\ref{fig:opt_eg}. As a result the dual function is non-differentiable and the bisection search may not converge~\cite{farias_2011}. (ii)~If $\lambda$ belongs to the sub-differential at a point $b_i$ of non-differentiability of either $f_i$ or $h_i$, the values of the gradient of $f_i\circ h_i$ may be arbitrarily different at $(b_i + \epsilon)$ and $(b_i-\epsilon)$ for an arbitrarily small $\epsilon$. This can also be seen in Fig~\ref{fig:opt_eg}. We use Algorithm~1 to compute the optimal solution of problem~(\ref{opt:cvx}). The convergence analysis is almost identical to that in Sec. 6 in~\cite{farias_2011}. An accurate solution can typically be computed in about 10 iterations.
\IncMargin{1em}
\begin{algorithm}
\DontPrintSemicolon
{\bf Given} starting value of $\underline{\lambda}$, $\overline{\lambda}$, $\underline{b}$, $\overline{b}$ and tolerance $\epsilon$. \;
\Repeat{$|\underline{\lambda}-\overline{\lambda}|<\epsilon$}{
\emph{Bisect:} $ \lambda = (\underline{\lambda} + \overline{\lambda})/2$.\;
\emph{Allocate bandwidth for all $i$:}\;
\eIf{$\lambda > \max\partial f_i(0)\max \partial h_i(0)$}{
set $b_i = 0$.\;
}{$b_i$ is such that\;
  \begin{equation}\begin{aligned}
  \label{eqn:bw_alloc}
  \smash{\lambda \in} &\smash{\left[\min \partial f_i(r_i) \times \min\partial h_i(b_i),\right.} \\
  &\smash{~~\left. \max \partial f_i(r_i) \times \max\partial h_i(b_i)  \right]}
  \end{aligned}\end{equation}\;
  \vspace{-0.6cm}
  where $$\smash{r_i = \left(b_i\psi\left(\frac{G_i(t) \min(\gamma_i b_i, P) }{Ib_i}\right) \right)}$$\;
\vspace{-0.6cm}
}

\emph{Update:} if $\sum_{i=1}^N b_i - B > 0$, $\underline{\lambda} = \lambda$, $\underline{b} = b$, else
$\overline{\lambda} = \lambda$, $\overline{b}=b$.\;
\BlankLine
}
\BlankLine
{\bf Feasible Solution:}
\eIf{$\sum_{i}\underline{b_i} - \sum_{i}\overline{b_i} > 0 $}{
set $\alpha = \frac{B - \sum_{i}\overline{b_i}}{\sum_{i}\underline{b_i} - \sum_{i}\overline{b_i}}$.\;
}{set $\alpha = 0$.
}
$ b = \alpha \underline{b} + (1-\alpha)\overline{b}$
\caption{Bisection search for optimal $\lambda$}\label{alg:bw_alloc}
\end{algorithm}
\DecMargin{1em}

The starting values of $\overline{\lambda}$ and $\underline{\lambda}$ can
be generated using the following simple lemma (proof is straightforward and omitted); the values of $\overline{b}$ and $\underline{b}$ are obtained by repeating the
\emph{Allocate Bandwidth} step in Algorithm~1 for dual variables $\overline{\lambda}$ and $\underline{\lambda}$, respectively.
\begin{lemma}
The optimal dual variable $\lambda^\star$ satisfies\\ $\underline{\lambda}\leq \lambda^\star\leq\overline{\lambda}$ where
\[
\begin{aligned}
\overline{\lambda} &= \max_{i=1,\hdots,N}\left[ \psi \left(\frac{G_i(t) \gamma_i}{I}\right)\max\partial f_i(0)\right]\\
&\\
\underline{\lambda} & = \left[\psi \left(\frac{G_i(t) P}{IB}\right) -\frac{G_i(t)P}{B}\psi'\left(\frac{G_i(t) P}{IB}\right)\right]\\
&\times \max\partial f_i\left(B\psi \left(\frac{G_i(t) P}{IB}\right) \right), ~~\text{for some } i
\end{aligned}
\]
\end{lemma}

The main computational step in each iteration of Algorithm~1 entails solving~(\ref{eqn:bw_alloc}) $N$ times -- we now show
this can be done in $O(\log L)$ time when the reward function $f_i$ for user $i$ is non-differentiable at at most $L$ points. The composition of function $f_i$ with $h_i$ is a concave function as shown in Lemma~(\ref{lem:cvx}).
Hence, to compute the bandwidth allocation for UE $i$
as given in equation~(\ref{eqn:bw_alloc}), we can use a bisection on $b_i$.
First we obtain
how many packets should be served fully such that the corresponding bandwidth required, $b_i$, satisfies equation~(\ref{eqn:bw_alloc}) in $O(\log L)$ time. Then, we compute $b_i$.

We compute the range of subgradients for packet $\eta$ as
\begin{equation}\begin{aligned}
\label{eqn:sg}
\underline{b} = h_i^{-1}\left(  \frac{\sum_{k=1}^{\eta-1} s_i}{\Delta}\right),\quad
\overline{b} = h_i^{-1}\left(  \frac{\sum_{k=1}^{\eta} s_i}{\Delta}\right)\\
SG(\eta) = d_i(\eta)[\min\partial h_i(\overline{b}), \min\partial h_i(\underline{b}) ]
\end{aligned}
\end{equation}
where we recall $d_i(\eta)$ and $s_i(\eta)$ are the delay and size for $\eta$th packet queued at UE $i$. Note that the inverse of $h_i$ is simple when $b_i < P/\gamma_i$; otherwise it can be computed via bisection.

\IncMargin{1em}
\begin{algorithm}
\DontPrintSemicolon

{\bf Initialize:} $\underline{\eta}=0$ $\overline{\eta}=\pi_i$, where we recall $\pi_i$ is the number of packets queued at UE $i$.\;
\Repeat{$\underline{\eta} = \overline{\eta}$}{
\emph{1. Bisect:} $ \eta = \left\lfloor (\underline{\eta} + \overline{\eta})/2\right \rfloor$.\;
\emph{2. Compute subgradient range $SG(\eta)$}\;
\emph{3. Update:} If $\min SG(\eta)>\lambda$, then $\underline{\eta} :=\eta$, else if
$\max SG(\eta)<\lambda$, then $\overline{\eta} :=\eta$, else $\overline{\eta}, \underline{\eta} :=\eta$.\;
}
\caption{Bisection for Number of Packets}\label{algo:invBisection}
\end{algorithm}
\DecMargin{1em}

The number of packets to be served completely is $\eta=\underline{\eta} -1$. Now we show how to compute
the bandwidth allocation $b_i$. Note that $h_i$ has at most one point of discontinuity, say $\hat{b}_i$.
If $\underline{b}\leq \hat{b}_i \leq \overline{b}$ for $\eta=\underline{\eta} -1$ in~(\ref{eqn:sg}), then
$b_i=\hat{b}$ if $\lambda/d_i(\eta)\in \partial h_i(\hat{b_i})$; else update
$\underline{b}$ or $\overline{b}$ appropriately.
\IncMargin{1em}
\begin{algorithm}
\DontPrintSemicolon

{\bf Given} tolerance $\mu$, $\underline{b}$, $\overline{b}$.\;
\Repeat{$|\underline{b} - \overline{b}|<\eta$}{
\emph{1. Bisect:} $ b =  (\underline{b} + \overline{b})/2$.\;
\emph{2. Update:} If $h_i'(b)>\lambda/d_i(\eta)$, then $\underline{b} :=b$, else if
$h_i'(b)<\lambda/d_i(\eta)$, then $\overline{b} :=b$.
}
\caption{Computation of RB assignment}\label{algo:rbAssign}
\end{algorithm}
\DecMargin{1em}

A similar method can be used for best effort traffic and the analysis is omitted here due to lack of space.

\section{Simulation Results}
\label{sec:results}

\subsection{Simulation Framework}

The algorithms in the previous section were simulated using a detailed system simulator where the MAC layer signalling was modeled faithfully, and the
PHY layer performance was abstracted via modeling of fading channels, transmission power, and capacity computations as in~\cite{modeling3},\cite{eval_tr}. A hexagonal regular cell layout with three sectors per site was simulated with the parameters as noted in Table~\ref{table:simParam}. For fractional power control parameter values ($P_0 = -60$ dBm, $\alpha = 0.6$) similar to those in~\cite{fpp_lte}, a 19 cell (57 sector) simulation with wrap around was first performed to determine the interference over thermal (IoT) at the base-station of a cell to be 6~dB on an average. In subsequent simulations, only one cell was simulated with the IoT assumed to be constant in time and frequency. This drastically reduces the simulation time while still accounting for the inter-cell interference.

The time varying channel gains, $G_i$'s, were assumed to be measured perfectly at the base-station in each subframe. The MCS was picked on the basis of the channel gain from the UE and a
rate adaptation algorithm to target an average of two H-ARQ transmissions for successful decoding was used. We use the mutual information effective SINR metric (MIESM)~\cite{BruAstSal05}; we first obtain the effective SINR according to the modulation alphabet size and then use that value to simulate an event of packet loss according to the packet error rate for the effective SINR. We model the timelines for Scheduling Request (SR), resource grants, Hybrid-ARQ, ACK/NACKs, and BSR as described in Sec.~\ref{sec:system}. We assume error free transmission of control messages in our simulations.

\begin{table}[!t]
\begin{tabular}{|l|l|}
\hline
Parameter & Value\\
\hline
Channel Profile & ITU-T PedA\\
\hline
Mobile Speed & 3 km/hr\\
\hline
Log-Normal Shadowing & $\sigma=$8.9 dBm\\
\hline
Intra-site Shadowing Correlation & 1.0\\
\hline
Inter-site Shadowing Correlation & 0.5\\
\hline
Cell Radius & 1 km\\
\hline
No. of UEs/cell  & 20\\
\hline
No. of RBs & 110\\
\hline
Max UE Tx Power & 23 dBm \\
\hline
No. of Tx \& Rx Antenna & 1 \\
\hline
eNB \& UE Antenna Gains & 0 dBi\\
\hline
Thermal Noise Density & -174 dBm/Hz\\
\hline
BSR periodicity & 5 ms\\
\hline
max. number of retransmission & 6\\
\hline
\end{tabular}
\caption{Simulation Parameters}
\label{table:simParam}
\end{table}

We focus on delay QoS traffic and consider two different models~\cite{traffic}. 
\emph{Live Video:} This is an ON-OFF Markov process with fixed packet size is used for the live video traffic model. The Markov process dwells in either state for 2 seconds and when in the ON state, generates a packet every 20 ms.
\emph{Streaming Video:}
\label{subsub:streaming_video}
Here, both the packet interarrival times and the packet sizes are independently drawn from truncated Pareto distributions. The number of arrivals in a frame length of 100 millisecond is fixed at 8, while their interarrival times are drawn from a truncated Pareto distribution with exponent 1.2 and truncation to [2.5 ms - 12.5 ms]. We use an exponent of 0.7 for the packet size distribution with varying values for the truncation limits, so as to control the mean data rate. For example, to get a mean rate of 500 kbps, we fix the limits to [215 bytes - 1500 bytes].


In order to map the optimal resource allocation computed using Algorithm~\ref{alg:bw_alloc} to actual RB grants we use a heuristic which ranks the users in decreasing order of marginal reward times spectral efficiency when very small amount of bandwidth is given to the user. The RB allocation is then done in the order of the rank, with each bandwidth amount (as per Algorithm~\ref{alg:bw_alloc}) mapped to an available segment of closest size.

\subsection{Results}
We consider two topologies for simulation: a \emph{macro-cell} with the path loss between the base station and UEs randomly selected between 100 dB and 135 dB~\cite{modeling3}, \emph{micro-cell} with path loss in the range 107 dB to 115 dB. We simulate three scheduling algorithms: (i)~\emph{Iterative Delay} which maximizes the reward function in Sec.~\ref{subsec:rew_delay}, (ii)~\emph{Iterative Queue} which minimizes sum-of-squares of queue lengths as in~\cite{sadiq_2009} and similar to~\cite{bodas_2011}, (iii)~\emph{non-iterative maximum weight} where a UE with the highest queue length times spectral efficiency for first RB is allocated bandwidth until the queue is drained or the UE becomes power limited before allocation to the next UE. We note that the computational algorithms in this paper are applicable to computing resource allocation for scheduling policies (i) and (ii), and that policies similar to (iii) do not consider the \emph{change} in reward function of the UE in a given subframe.

\noindent \subsubsection{\bf{Macro cell Topology}}
We consider 20 UEs with a mix of live video and streaming video traffic. Since live video has a tighter requirement for packet delays, we bias the scheduler to assign live video users 5x priority compared to streaming video users for same packet delay. Simulations were performed for low load and high load cases:\\
(1)~\emph{High Load:} 5 UEs have live video traffic, each with a mean rate of 300 kbps. For the other 15 UEs with streaming video traffic, we mimic an adaptive-rate streaming mechanism in which the data rate for each user depends on the quality of its channel to the base-station, i.e. a user close to the base-station transmits a better quality video compared to a cell-edge user. For simulating high-load, the truncation parameters mentioned in Section~\ref{subsub:streaming_video} are varied for each UE such that they generate traffic at 80\% of the average data rate they received with full buffer traffic.\\
(2)~\emph{Low Load:} 5 UEs have live video traffic with a mean rate of 200 kbps. The UEs with streaming video traffic are now set to operate at 40\% of their full buffer average data rate.

\begin{figure}[!t]
\par
\begin{center}
\includegraphics[width=0.45 \textwidth]{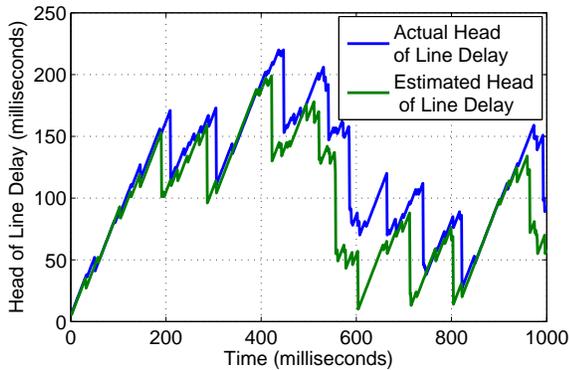}
\caption{HoL delay estimation performance}
\label{fig:qEstimation}
\end{center}
\end{figure}

We first study the performance of the delay estimation mechanism described in Section~\ref{subsec:inferring_arrivals}. Figure~\ref{fig:qEstimation} shows the estimated head of line (HoL) delay and the actual HoL delay at a UE over a period of 1 second. The estimated values can be seen to follow the actual delays but the accuracy is limited by the granularity of BSR messages, i.e., if there are multiple arriving packets between two successive BSR messages, the packets are bundled as one in our mechanism resulting in relatively small errors in HoL estimation.

\begin{figure}[!t]
\par
\begin{center}
\includegraphics[width=0.5\textwidth]{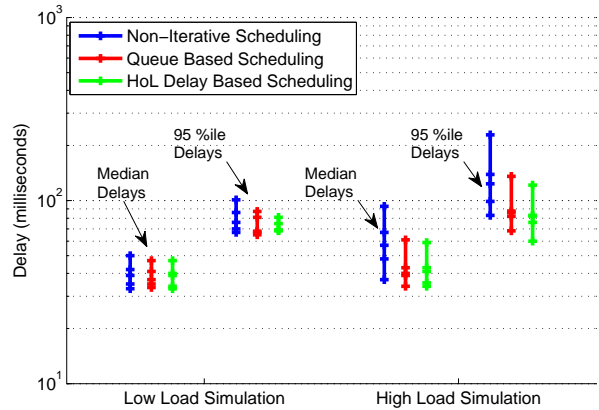}
\caption{Live video users: delay performance}
\label{fig:delayLive}
\end{center}
\end{figure}
\begin{figure}[!t]
\begin{center}
\includegraphics[width=0.5\textwidth]{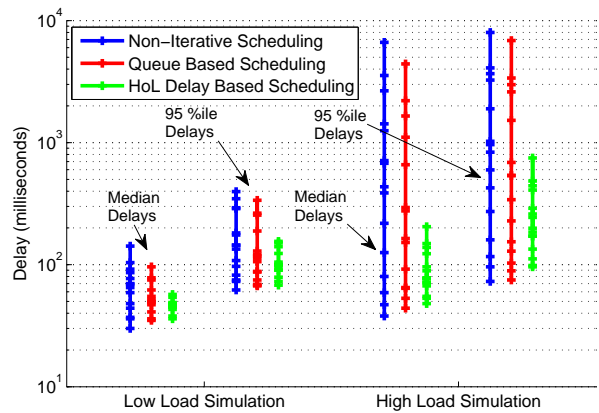}
\caption{Streaming video users: delay performance}
\label{fig:delayStream}
\end{center}
\end{figure}
\begin{figure}[!t]
\begin{center}
\includegraphics[width=0.5\textwidth]{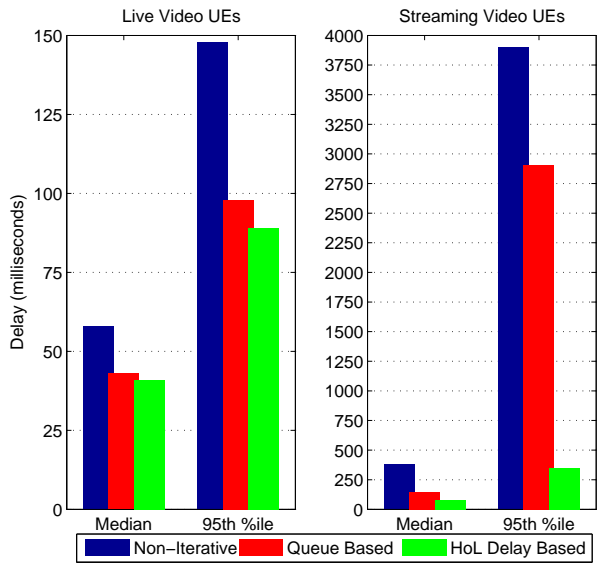}
\caption{Cell-wide delay performance of all packets in macro cell simulation}
\label{fig:allUEs}
\end{center}
\end{figure}

Next we show the performance of the head of line delay based scheduling scheme computed as the solution to the optimization problem in~(\ref{opt:cvx}) with the reward function in~(\ref{eqn:rew_qos}). Figure~\ref{fig:delayLive} shows the median and 95th percentile delays of the live video UEs for the two baseline and the head of line delay based schedulers for low and high loads. The delays experienced by the live video users are consistently less in the case of HoL delay based scheduling with the non-iterative scheme resulting in an average 95th percentile delay 1.6x higher than with the HoL delay scheduling. The queue based scheme also results in slightly higher delays, on an average 1.1x compared to 95th percentile delays for HoL scheduling. A more pronounced improvement is observed for the streaming video users, as shown in the delay plots in Figure~\ref{fig:delayStream}. In this case, the non-iterative and queue based schemes result in 6.2x and 5x more delays compared to HoL delay scheduling in terms of 95th percentile latencies. Finally, Figure~\ref{fig:allUEs} shows the combined delay numbers for uplink packets from all the UEs in the high load simulation. As can be seen from the figure, the iterative queue based and delay based schemes result in similar delays for live video users due to preferential assignment. However this results in large delays for the streaming video users for both non-iterative and queue based schemes: close to 11x and 8x respectively compared to HoL delay based scheduling in terms of 95th percentile delays. Thus, leveraging the approximate packet delays obtained via our method leads to significant performance improvement over queue based scheduling. Moreover, even for the queue based scheduler, the computational methods in this paper are very useful.

\noindent \subsubsection{\bf{Micro cell Topology}}
\begin{figure}[!t]
\begin{center}
\includegraphics[width=0.5\textwidth]{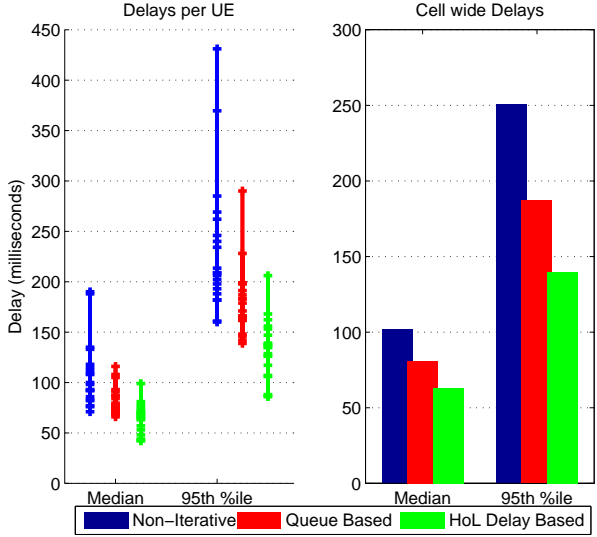}
\caption{Individual and Cell-wide Delay performance for micro cell simulation}
\label{fig:microCell}
\end{center}
\end{figure}

In order to compare these scheduling schemes in a smaller cell topology, we ran a second simulation with 20 UEs located within a region with path loss 107-115 dB from the base station. Each UE, in this simulation, carries streaming video traffic with the mean data rate randomly selected between 300-2000 Kbits/sec. Decoupling the mean traffic rate with the path loss highlights the relative  performance of the scheduling algorithms in real deployments where prior knowledge of user demand is rarely known. Individual and cell wide delay numbers are shown in Figure~\ref{fig:microCell}, which shows that 95th percentile delays for non-iterative and queue based schemes are 1.8x and 1.4x more than those for the HoL delay based scheduling.

\section{Frequency Selective Resource Allocation}
We extend the analysis in~\cite{madan_2011} for frequency selective fading to concave functions $f_i$ (such as the delay based reward function) which are thrice continuously differentiable everywhere except at $L$ points where they are only continuous. We can re-write such a function as
\[ \smash{f_i(r_i) = \sum_{l=1}^L f_{il}\left(\min\left(\rho_l-\rho_{l-1}, \left[r_i - \rho_{l-1}\right]_+\right)\right)}\]
where $0\leq \rho_1 < \hdots <\rho_L$ are the points of non-differentiability and $f_{il}:\mathbb{R}_+\mapsto \mathbb{R}$ are thrice continuously differentiable concave functions defined as
\[\smash{f_{il}(x) = f_i(\rho_{l-1}+x) - f_i(\rho_{l-1}), ~x\in[0,\rho_l-\rho_{l-1}],l\geq 1}, \rho_0=0,\]
and satisfy
 $$
\smash{f_{il}'(x) < f_{i,l-1}'(y), \quad l>1,x\in[0,\rho_l-\rho_{l-1}],~~y\in[0,\rho_{l-1}-\rho_{l-2}].}
 $$
 We also assume $x\psi^{-1}(y/x)$ is concave for all $(x,y)>0$; this is true for example, when $\psi$ is the Shannon capacity formula, and for practical
M-QAM schemes.

Consider the following convex optimization problem over $\tilde{r}_{il}$'s, $r_{ij}$'s (rate for user $i$ on sub-band $j$), and $b_{ij}$'s (bandwidth for user $i$ on sub-band $j$):
\begin{equation}
\label{opt:enb_comb}
\begin{aligned}
\text{max.}~~ & \sum_{i=1}^N\sum_{l=1}^L f_{il}(\tilde{r}_{il}), \\
\text{s.t.}~~ & \sum_{l=1}^L \tilde{r}_{il} \leq \sum_{j=1}^M r_{ij},~~\forall i, ~~ \tilde{r}_{il}\leq \rho_{l}-\rho_{l-1},~\forall i,l \\
                & \sum_{i=1}^Nb_{ij} = B, ~~ \forall j,\\
                 & \sum_{j=1}^M \frac{b_{ij}(N_0 + I_j)}{G_{ij}}\psi^{-1}\left(r_{ij}/b_{ij} -
                 1\right)\leq P, ~~\forall i,\\
                 & r_{ij}\leq b_{ij}\psi\left(\frac{G_{ij}\gamma_{ij}}{N_0 + I_j}\right),~~ r_{ij}, b_{ij}\geq 0,
~\forall i,j. \\
\end{aligned}
\end{equation}
The first constraint implies that the total rate for a user is the sum of rates over sub-bands, the second constraint is on total bandwidth allocation in a sub-band, third constraint is on peak power at the UE in a subframe, and the fourth constraint models fractional power control.
The following lemma follows easily from the construction of the $f_{il}$s:
\begin{lemma}
If $(r^\star_{ij}, b^\star_{ij})$ is a solution to the optimization problem~(\ref{opt:enb_comb}), then $\sum_{i}f_i\left(\sum_{j}r_{ij}^\star\right)$
is the maximum sum reward for any feasible resource allocation.
\end{lemma}

General purpose interior points methods to solve the above optimization problem have a complexity of $O(NM+NL)^3$ per iteration -- we exploit the structure to reduce it to $O(N(L^2+M^2))$. Note that in practice $L$ and $M$ are much smaller than $N$. In order to construct a solution for which the bandwidth allocation is contiguous in frequency to satisfy the SC-FDMA requirements, we can use the heuristic in~\cite{madan_2011}. The main computation to solve~(\ref{opt:enb_comb}) is to determine the Newton step at each iteration which entails solving a set of linear equations of the form (we omit the details due to lack of space, the exact expressions can be obtained following the steps in~\cite{boyd_2004}):
\[
\left[
\begin{array}
{cccccc}
      H_{1}  & &  & & \vline & A^T\\
      & & \ddots  & & \vline&  \\
      & & & H_{N}  & \vline& \\
      \hline &&&&\vline & \\
      && A &&\vline& 0
\end{array}
\right]
\left[\begin{array}{c}
x_1\\
\vdots\\
x_N\\
y
\end{array}\right]
=
\left[\begin{array}{c}
a\\
b
\end{array}\right]
\]
where $H_i\in \mathbb{R}^{(L+M)\times(L+M)}$, $A\in \mathbb{R}^{M\times N(L+M)}$, $x_i\in \mathbb{R}^{L+M}$, $a\in \mathbb{R}^{N(L+M)}$, $y,b\in \mathbb{R}^{M}$. We first eliminate the $x_i$'s as
\[
x_i = H_i^{-1}\left(a - A^T_{(L+M)(i-1)+1:(L+M)(i)} y\right)
\]
where $A^T_{k:m}$ is the submatrix of $A^T$ given by rows $k$ to $m$. We invert $H_i^{-1}$ in $O(L^2+M^2)$ time, solve for $y$ in $O(M^3)$ time ($M$ linear equations in $M$ variables), and back-substitute $y$ to obtain $x$. To invert $H_i$, we note that it decomposes as

\[\begin{aligned}
H_i&=\left[
\begin{array}
{ccccccccc}
      K_{1}  &&&&&&&&\\
     & & \ddots  & & &&&&  \\
     & & & K_{M}  & &&&& \\
      &&&& g_{1}  & &  & & \\
      &&&& & & \ddots  & &  \\
      &&&& & & & g_{L}  &  \\
\end{array}
\right ]
+ g_ig_i^T  \\
&\quad+\left[
\begin{array}
{cc}
h_ih_i^T &0\\
0&0
\end{array}
\right ] +
\left[
\begin{array}
{cc}
0 &0\\
0& c_ic_i^T
\end{array}
\right ]
\end{aligned}
\]
where $g_i\in\mathbb{R}^{L+M}$, $h_i\in\mathbb{R}^{L}$, $c_i\in\mathbb{R}^{M}$. Using the matrix inversion lemma we can invert $H_i$ in $O(L^2+M^2)$ time.

\section{Conclusions}
We designed a general computational framework in this paper to enable a wide array of online scheduling policies in a computationally efficient manner. We modeled the constraints due to fractional power control, and formulated an optimization problem with non-differentiable objective function. We showed how to estimate the packet delays on the uplink via the BSR reports, and proposed a novel scheduling policy based on packet delays. Numerical results demonstrated that using packet delay estimates for the uplink can lead to significant reduction in packet delays as compared with a queue length based scheduler. There are many interesting directions for future work. For example, we can further study the connections with the work in~\cite{sharma_2011},\cite{bodas_2011}. In terms of implementation, an interesting question is whether we can design approximation algorithms for the uplink bandwidth packing problem which are optimal according to some metric.

\footnotesize
\bibliographystyle{IEEEtran}
\bibliography{fast_tdma}

\end{document}